\documentclass[12pt,a4paper,twoside]{amsart}
\usepackage[english]{babel}
\usepackage{amsmath}
\usepackage{amssymb}
\usepackage[latin1]{inputenc}  %<--- fuer UNIXe
\theoremstyle{plain}
\newtheorem{theorem}{Proposition}[section]
\newtheorem{thm}[theorem]{Theorem}
\newtheorem{lemma}[theorem]{Lemma}

\newtheorem{corol}[theorem]{Corollary}
\renewcommand{\H}{\mathcal H}
\newcommand{\N}{\mathbb N}

\newcommand{\C}{\mathbb C}
\newcommand{\A}{\mathcal A}
\newcommand{\D}{\mathcal D}
\newcommand{\F}{\mathcal F}
\newcommand{\B}{\mathcal B}
\newcommand{\la}{\big\langle}
\newcommand{\ra}{\big\rangle}
\renewcommand{\L}{\mathcal L^+(\D)}
\begin{document}

\title[]{Zero Product Preservers and Orthogonality Preservers in Algebras of Unbounded Operators}

\author[]{WERNER TIMMERMANN}

\address{Institut f\"ur Analysis\\
         Technische Universit\"at Dresden\\
         D-01062 Dresden, Germany}
\email{timmerma@math.tu-dresden.de}
\begin{abstract}

Applying a result of abstract ring theory we get that bijective additive mappings on standard algebras of unbounded operators preserving zero products are multiples of ring isomorphisms. \\
The structure of additive bijective mappings on certain classes of standard algebras of unbounded operators preserving orthogonality in both directions is also  investigated. The results are quite similar to those for algebras of bounded operators.
\end{abstract}

\thanks{\hspace{-0.5cm}Keywords: algebras of unbounded operators, preservers }

\maketitle
\section{Introduction}
Linear preserver problems concern the characterization of linear maps between algebras that roughly speaking preserve certain properties of some elements of the algebras. Such problems were studied in matrix theory during the last century starting with the paper of Frobenius [4].\\
In the last decades there can be observed a growing interest in similar questions on abstract algebras or rings and on operator algebras over infinite dimensional spaces. These investigations led in a natural way to the study of mappings which are not linear but are merely additive. One of the first paper in this spirit is the classical work of Kaplansky [6]. One of the striking features is the application of results of abstract ring theory to problems on algebras of bounded operators on Banach or Hilbert spaces. \\
It turned out that several of these results are valid with necessary modifications also for algebras of unbounded operators on Hilbert spaces (see for example [9 -- 12]). \\
The present paper is exactly in this spirit. It deals with two special preserver problems: preservers of zero products and preservers of orthogonality. The results are based on ring theoretical theorems concerning these topics.\\
We fix the necessary notions and notation. Let  $\A, \B$ be rings. A mapping $\Phi:~ \A \to \B$ is said to preserve zero products if $XY = 0 ~ (X,Y \in \A)$ implies $\Phi(X)\Phi(Y) = 0$. Let $\A, \B$ be *-algebras. A mapping $\Phi: \A \to \B$ preserves orthogonality if $A^*B = AB^* = 0 ~(A,B \in \A)$ implies $\Phi(A)^*\Phi(B) = \Phi(A)\Phi(B)^* = 0$. If $\Phi$ is bijective then $\Phi$ is said to preserve orthogonality on both direction if
$$A^*B = AB^* = 0 \quad \Longleftrightarrow \quad \Phi(A)^*\Phi(B) = \Phi(A)\Phi(B)^* = 0.$$
Further we need some notions on algebras of unbounded operators.
 A standard reference  is [8].\\ 
Let $\mathcal D$ be a dense linear manifold in a Hilbert space $\mathcal H$ with scalar
product $\la ~,~ \ra $ (which is supposed to be conjugate linear in the first and linear in the second component). The set of linear operators $\L = \{ A: A \D \subset \D, A^* \D \subset \D\}$ is a $\ast$-algebra with respect to the natural operations and the involution $A \rightarrow A^+ = A^*|\D$. The graph topology $t$ on $\D$ induced by $\L$ is generated by the directed family of seminorms $\phi \rightarrow ||\phi ||_A = || A \phi ||,~ \forall ~  A \in \L, ~  \phi \in \D$. $\D$ is called an (F)-domain, if $(\D, t)$ is an (F)-space.\\
 Remark that in this case the graph topology $t$ can be given by a system of seminorms $\{ \| \cdot \|_n = \|A_n \cdot \|, n \in  \N, A_n \in \L\}$ with:
 \begin{equation}
A_1 = I, \  A_n = A^+_n, \   \|A_n \phi\| \leq \|A_{n+1} \phi \| \quad \text{for all} \ \phi \in \D, n \in  \N. 
\end{equation} 
 A standard (*-) operator algebra (on $\D$) is a (*-) subalgebra $\mathcal A(\D) \subset \L$ containing the ideal $\mathcal F(\D) \subset \L$ of all finite rank operators on $\D$. Note that every rank-one operator $F \in \F(\D)$ has the form $F = \psi \otimes \phi, ~ \phi, \psi \in \D$, where $F(\chi) = \la \phi, \chi\ra \psi$.\\
Every standard operator algebra $\A(\D)$  is prime. Remember that an algebra (or a ring)  $\A$ is prime if $X\A Y = 0$ implies $X = 0$ or $Y = 0$.\\
The paper is organized as follows. In section 2 we deal with additive bijective mappings between standard operator algebras preserving zero products. The corresponding result is a simple application of an abstract ring theoretical result [2]. It appears that such mappings are scalar multiples of ring isomorphisms.\\
Section 3 is devoted to additive bijective mappings preserving orthogonality. The structure of such mappings is clarified. The results are quite similar to those for algebras of bounded operators given in [5].

\section{Mappings preserving zero products}
For the extensive literature on this topic see for example [2,3] and the references therein. In this section we prove the following theorem.
\begin{thm}
Let $\D \subset \H$ be an (F)-domain and let $\A, \B$ be standard *- operator algebras on $\D$. If $\Phi:~ \A \to \B$ is an additive, bijective mapping that preserves zero products, then there are $c \in \C$ and $T: \D \to \D$ bijective and either linear or conjugate linear such that
$$\Phi(A) = c TAT^{-1} \qquad (A \in \A).$$
If $T$ is linear, then $T \in \L$
\end{thm}
The main part of the proof is contained in the following  theorem from ring theory (Theorem 1 in [2]).\\

{\bf Theorem A} \emph{Let $\A$ and $\B$ be prime rings and $\Phi: \A \to \B$ a bijective additive mapping such that $\Phi(A) \Phi(B) = 0$ for all $A,B \in \A$ with $AB = 0$. Suppose that the maximal right quotient ring $Q(\A)$ of $\A$ contains a nontrivial idempotent $E$ such that $E\A \cup \A E \subset \A$.\\
i) If $1 \in \A$, then $\Phi(AB) = \lambda \Phi(A)\Phi(B) $ for all $A,B \in \A$ where $\lambda = \frac{1}{\Phi(1)} \in Z(\B)$ the center of $\B$. In particular, if $\Phi(1) = 1$, the $\Phi$ is a ring isomorphism from $\A$ onto $\B$.\\
ii)  If deg($\B) \geq 3$, then there exists $\lambda \in C(\B)$, the extended centroid of $\B$, such that $\Phi(AB) = \lambda \Phi(A)\Phi(B)$ for all $A,B \in \A$}\\

The definitions and basic properties of the maximal quotient ring and the extended centroid can be found in [1]. A prime ring $\A$ is called centrally closed if $C(\A)$ is trivial. We need the following characterization of centrally closed prime algebras to prove that every standard operator algebra on $\D$ is centrally closed (thanks to M. Bre\u{s}ar for this information). Let $\A$ be a prime algebra over $\C$. Then $\A$ is centrally closed if and only if the following holds:\\
If $\mathcal I \subset \A$ is a nonzero ideal and if there is an additive mapping $F: \mathcal I \to \A$ such that
\begin{equation} 
F(UX) = F(U)X \text{ and } F(XU) = XF(U)
\end{equation}
for all $U \in \mathcal I, X \in \A$, then there is a $\lambda \in \C$ such that $F(U) = \lambda U$ for all $U \in \mathcal I$.\\ In the next lemma we use a modification of the notion of a double centralizer. A pair of (linear or additive) mappings $L,R: \mathcal I \to \A$ is called a double centralizer if
$$ L(XY) = L(X)Y, \quad R(XY) = XR(Y) \quad \text{ and } \quad XL(Y) = R(X)Y  \text{ for all } X,Y \in \A$$
The structure of double centralizers on standard operator algebras on $\D$ was described in [10] Proposition 3.4 as follows: there is a $T \in \L$ such that $L(A) = TA,~  R(A) = AT$. The same proof is valid also for the modified situation descibed above. We use this result to prove the following lemma.
\begin{lemma}
Every standard operator algebra  on $\D$ is centrally closed.
\end{lemma}
\begin{proof}{}
Let $\mathcal I \subset \A$ be an ideal and $F: \mathcal I \to \A$ an additive mapping such that (2) is satisfied. Note that $\F(\D) \subset \mathcal I \subset \A$. The pair $(F,F)$ is a double centralizer in the sense described above. Consequently there is a $T \in \L$ such that $F(A) = TA, F(A) = AT$ for all $A \in \mathcal I$. So, $T$ commutes with all operators from $\mathcal I$, in particular, with all operators from $\F(\D)$. But this implies $T = \lambda I$ for some $\lambda \in \C$.
\end{proof}

{\bf Proof of Theorem 2.1 :} Apply Theorem A, ii) to get $\Phi = \frac{1}{\lambda}\Psi$ with a ring isomorphism $\Psi$. The structure of ring isomorphisms between standard *-operator algebras on (F)-domains follows directly from Theorem 3.1 in  \cite{ti03a}. Namely, there exists a bijective either linear or conjugate linear $T: \D \to \D$ such that $\Psi(A) = TAT^{-1} ~ (A \in \A)$. This concludes the proof.\\

\section{Mappings preserving orthogonality}
In this section we describe the structure of orthogonality preserving mappings on several standard operator algebras. As a corollary we obtain an unbounded version of a result of L. Moln\'ar [7].
\begin{thm}
Let $\D \subset \H$ be an (F)-domain and let $\A \subset \L$ be one of the following standard operator algebras:\\
a) $\A = \F(\D)$,\\
b) $\A$ is a unital standard *- operator algebra,\\
c) $\A \not= \F(\D)$, $\A$ a *-ideal of $\L$.\\
Assume that $\Phi: \A \to \A$ is an additive bijection preserving orthogonality in both directions. Then $\Phi$ has one of the following forms:\\
i) There exist a nonzero constant $c$ and operators $U,V: \D \to \D$, both either unitary or antiunitary such that
\begin{equation}
 \Phi(T) = cUTV \qquad (T \in \A)
\end{equation}
or\\
ii) There exist a  nonzero constant $c$ and operators  $U,V: \D \to \D$, both either unitary or antiunitary such that
\begin{equation}
 \Phi(T) = cUT^+V \qquad (T \in \A)
\end{equation}
\end{thm}
\begin{proof}{}
The main step consists in proving that $\Phi$ preserves rank-one operators in both directions. For this we treat cases a) - c) separately.\\
Case a): This can be done as in [5] using the polar decomposition of $A \in \F(\D)$. Remark that in the unbounded case the situation is more complicated.  For $\A \not= \F(\D)$ one can not be sure that the polar decomposition can be performed within $\A$ (even not within $\L$).\\
Case b): Let $F = \psi \otimes \phi$ be a rank-one operator and let $\H_0 := lin\{\phi, \psi\}$. Now we define an operator $Q \in \A$, which has corank one as follows.\\
 If dim $\H_0 = 1$ put
$$Q \chi = \begin{cases}  \chi & \quad \text{ for all } \chi \in \H_0^\perp \cap \D\\
0 & \quad  \text{ for all } \chi \in \H_0. \end{cases} $$   
If dim $\H_0 = 2$ put
$$Q \chi = \begin{cases}  \chi & \quad \text{ for all } \chi \in \H_0^\perp \cap \D\\
\la \phi_1, \chi\ra\psi_1 & \quad  \text{ for all } \chi \in \H_0. \end{cases} $$ 
where $\phi_1, \psi_1 \in \H_0$ such that $\phi_1 \perp \phi, \psi_1 \perp \psi$.  To get $Q  \in \A$ it is used that $I \in \A$.\\
Then $F^+Q = FQ^+ = 0$ and therefore $\Phi(F)^+\Phi(Q) = \Phi(F)\Phi(Q)^+ = 0$.\\ That means ran$\Phi(F) \perp$ ran$\Phi(Q)$ and ran$\Phi(F)^+ \perp $ ran$\Phi(Q)^+$.\\
Suppose that $\Phi(F)$ has rank larger than one. Then also $\Phi(F)^+$ has rank larger than one. Let $\rho_1, \rho_2 \in \text{ran} \Phi(F), \rho_1 \perp \rho_2$ and let $\chi_1, \chi_2 \in \text{ran} \Phi(F)^+, \chi_1 \perp \chi_2$. Put $S_i := \rho_i \otimes \chi_i$. Then $S_i: \D \to \text{ran} \Phi(F), S_i^+: \D \to \text{ran} \Phi(F)^+$ and $S_1, S_2$ are orthogonal. Obviously  $S_i^+\Phi(Q) = S_i\Phi(Q)^+ = 0$.\\
Now let $T_i \in \A$ such that $\Phi(T_i) = S_i$ Then $T_1T_2^+ = T_1^+T_2 = 0$ and $T_i^+Q = T_iQ^+ = 0$, i.e. ran $T_1 \perp$ ran $T_2$ and ran $T_i \perp$ ran $Q$. But this is a contradiction, because  ran $Q$ has codimension one.\\
Case c): Here we can argue similar to case b).\\
 Let $F = \psi \otimes \phi \in \F(\D)$, $\H_0 = $lin$\{\phi, \psi\}, \D_0 = \H_0^\perp \cap \D.$\\
Fix an operator $A \in \A$ with infinite-dimensional range. Let $(\phi_n) \subset $ ran$A$ be an orthonormal basis from $\overline{\text{ran} A}$ and let $(\psi_n) \subset \D_0$ be an orthonormal basis from $\H_0^\perp$. Now we use property (1) of the operators defining the topology $t$ in $\D$. One can find a positive, bounded  sequence $(\alpha_n)$ such that the operator $U$ defined by $U = 0$ on (ran $A)^\perp$ and $U \phi_n = \alpha_n \psi_n$ (and linear, bounded extension to $\overline{\text{ran}A}$) has the following properties. $U \H \subset \D, U^*\H \subset \D$, so $U\D \subset \D, U^*\D \subset \D$ and consequently the restriction of $U$ to $\D$ (also denoted by $U$) belongs to $\L$. Moreover by construction, $U \D \subset \D_0$ and ran $U$ is dense in $\H_0^\perp$. This implies that $UA \in \L$ and ran$UA$ is dense in $\H_0^\perp$. Analogously to case b) we define an operator $Q$ as follows.\\
If dim $\H_0 = 1$ put
$$Q \chi = \begin{cases} UA \chi & \quad \text{ for all } \chi \in \H_0^\perp \cap \D\\
0 & \quad  \text{ for all } \chi \in \H_0. \end{cases} $$   
If dim $\H_0 = 2$ put
$$Q \chi = \begin{cases}  UA\chi & \quad \text{ for all } \chi \in \H_0^\perp \cap \D\\
\la \phi_1, \chi\ra\psi_1 & \quad  \text{ for all } \chi \in \H_0. \end{cases} $$ 
where $\phi_1, \psi_1 \in \H_0$ such that $\phi_1 \perp \phi, \psi_1 \perp \psi$. \\
As in case b) this operator $Q$ can be used to prove that  $\Phi(F)$ is a rank-one operator.\\
So, in any of the cases a) - c) $\Phi$ restricted to $\F(\D)$ is a bijective additive mapping which preserves rank-one operators in both directions.  According to [12] we get the following structure of $\Phi$ on $\F(\D)$. Either there exist  a ring automorphism $h: \C \to \C$ and bijective additive mappings $B,C: \D \to \D$ with the properties $B(\lambda \phi) = h(\lambda)B\phi, C(\lambda \phi) = h(\lambda)C\phi$ such that
$$\Phi(\psi \otimes \phi) = B\psi \otimes C\phi \qquad (\phi, \psi \in \D)$$
or there exist a ring automorphism $k: \C \to \C$ and bijective additive mappings $B,C: \D \to \D$ with the properties $B(\lambda \phi) = k(\lambda)B\phi, C(\lambda \phi) = k(\lambda)C\phi$ such that
$$  \Phi(\psi \otimes \phi) = C\phi \otimes B\psi\qquad (\phi, \psi \in \D).$$
Now the proof can be completed as in the last part of the proof of Theorem 4 in  [5].
\end{proof}

As a corollary we get an unbounded version of Theorem 2 in [7]
\begin{corol}
Let $\D \subset \H$ be an (F)-domain and let $\A(\D)$ be a unital standard operator *-algebra. Suppose that $\Phi: \A \to \A$ is an additive  bijection such that\\
i) $A^+B = 0$ if and only if $\Phi(A)^+\Phi(B) = 0$ and \\
ii) $AB^+ = 0$ if and only if $\Phi(A)\Phi(B)^+ = 0$.\\
Then there exist a nonzero constant $c$ and operators $U,V: \D \to \D$, both either unitary or antiunitary such that 
\begin{equation}
\Phi(T) = c UTV \qquad (T \in \A(\D)).
\end{equation}
If in addition $\Phi(I) = I$ then there is a unitary or antiunitary operator $U: \D \to \D$ such that
\begin{equation}
\Phi(T) = UTU^+  \qquad (T \in \A(\D))
\end{equation}

\end{corol}
\begin{proof}{}
If $\Phi$ satisfies i) and ii) then $\Phi$ preserves orthogonality in both directions. Consequently $\Phi$ has one of the representations (3),(4) of Theorem 3.1.\\
 First we exclude representation (4). 
Suppose $A^+B = 0$. Then 
$$\Phi(A)^+\Phi(B) = \bar cV^+AU^+cUB^+V = |c|^2 V^+AB^+V.$$
But this expression is not necessarily zero because $AB^+$ must not be zero if $A^+B = 0$. So it remains the form (3).\\
Now let $\Phi(I) = I$. Then $\Phi(I) = cUV = I$. So $V = c^{-1} U^{-1} = c^{-1}U^+$ Therefore
$$\Phi(T) = UTU^+.$$
\end{proof}

\end{document}